\documentclass{article}

\usepackage{amsmath}

\usepackage[amsmath,thmmarks,hyperref]{ntheorem}

\usepackage{hyperref}
\hypersetup{
    colorlinks=true,
    linkcolor=blue,
    citecolor=blue,
    breaklinks
}

\usepackage{graphics,graphicx}

\usepackage{algorithm}

\usepackage[capitalize,nameinlink]{cleveref}

\usepackage{lipsum}
\usepackage{amsfonts}
\usepackage{graphicx}
\usepackage{epstopdf}
\usepackage{algpseudocode}
\DeclareGraphicsExtensions{.eps,.pdf,.png,.jpg}
\usepackage{subcaption}
\usepackage{booktabs}

\floatname{algorithm}{Algorithm}

\usepackage{xcolor, bm}

\newcommand{\scrF}{\mathcal{F}} % shortcut for script F
\newcommand{\reals}{\mathbb{R}} % shortcut for Real numbers
\renewcommand{\vec}[1]{\bm{#1}} % redefine vector notation to be bold
\DeclareMathOperator*{\argmin}{arg\,min} % declare \argmin like \min

\usepackage{enumitem}

\title{Data fusion for a multi-scale model of a wheat leaf surface: a unifying approach using a radial basis function partition of unity method%
\thanks{This work was funded by the Australian Research Council through the ARC Linkage Project LP160100707, and the associated industry partners Syngenta and Nufarm. Funding was also provided by the Australian Research Council through the ARC Training Centre for Multiscale 3D Imaging, Modelling and Manufacturing (M3D Innovation, project IC 180100008)}}

\author{Riley M. Whebell\thanks{School of Mathematical Sciences, Queensland University of Technology, Brisbane 4000, Australia}
\and Timothy J. Moroney\footnotemark[2]
\and Ian W. Turner\footnotemark[2]
\and Ravindra Pethiyagoda\thanks{School of Information and Physical Sciences, University of Newcastle, New South Wales 2308, Australia}
\and Marie-Luise Wille\thanks{ARC ITTC for Multiscale 3D Imaging, Modelling and Manufacturing, Queensland University of Technology, Brisbane, Australia; and Centre for Biomedical Technologies, School of Mechanical, Medical, and Process Engineering, Faculty of Engineering, Queensland University of Technology, 60 Musk Avenue, Kelvin Grove, Queensland, 4059, Australia}
\and Justin J. Cooper-White\thanks{School of Chemical Engineering, and Australian Institute for Bioengineering and Nanotechnology, The University of Queensland, Australia}
\and Arvind Kumar\footnotemark[5]
\and Philip Taylor\thanks{Syngenta, Jealott’s Hill International Research Centre, Bracknell, Berkshire RG42 6EY, U.K}
\and Scott W. McCue\footnotemark[2]
}

\usepackage{amsopn}

\hypersetup{
  pdftitle={Multi-scale model of a wheat leaf surface},
  pdfauthor={R. M. Whebell, et al.}
}

\begin{document}

\maketitle

\begin{abstract}
  Realistic digital models of plant leaves are crucial to fluid dynamics simulations of droplets for optimising agrochemical spray technologies. 
The presence and nature of small features (on the order of $100\mathrm{\mu m}$) such as ridges and hairs on the surface have been shown to significantly affect the droplet evaporation, and thus the leaf’s potential uptake of active ingredients. 
We show that these microstructures can be captured by implicit radial basis function partition of unity (RBFPU) surface reconstructions from micro-CT scan datasets. 
However, scanning a whole leaf ($20\mathrm{cm}^2$) at micron resolutions is infeasible due to both extremely large data storage requirements and scanner time constraints. 
Instead, we micro-CT scan only a small segment of a wheat leaf ($4\mathrm{mm}^2$). 
We fit a RBFPU implicit surface to this segment, and an explicit RBFPU surface to a lower resolution laser scan of the whole leaf.
Parameterising the leaf using a locally orthogonal coordinate system, we then replicate the now resolved microstructure many times across a larger, coarser, representation of the leaf surface that captures important macroscale features, such as its size, shape, and orientation. 
The edge of one segment of the microstruture model is blended into its neighbour naturally by the partition of unity method.
The result is one implicit surface reconstruction that captures the wheat leaf’s features at both the micro- and macro-scales.
\end{abstract}
\vspace{1em}
\\[1em]

\section{Introduction}
\label{sec:intro}
In agriculture, understanding the behaviour of water droplets on leaves is critical for improving droplet spray efficacy, for example in the application of pesticides.
Realistic virtual plant leaf surfaces are a key component of droplet spray efficacy models \cite{dorr_spray_2016, dorr_towards_2014, dorr_impaction_2015, mayo_simulating_2015, tredenick_evaporating_2021}.
It is known that droplet behaviour on leaves such as wheat is highly dependent on the micro-structures such as the ridges and hairs of the leaf \cite{tredenick_evaporating_2021}.
Hence for maximum realism, virtual leaf surfaces used in droplet spray models should be reconstructed from scanned data of real leaves, in order to faithfully represent the morphology of the plant.

\begin{figure}
    \centering
    \begin{subfigure}{0.3\linewidth}
        \includegraphics[width=\linewidth]{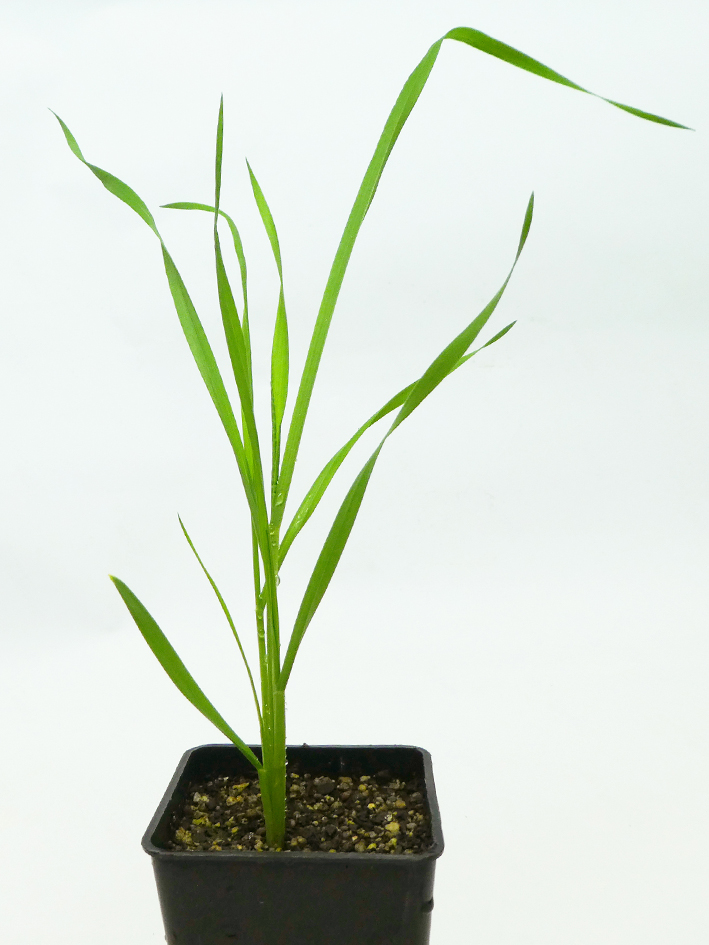}
        \caption{A photograph of a wheat plant, similar to that which was laser scanned for our macro-scale dataset.}
        \label{figs:wheatPhoto}
    \end{subfigure}
    \hspace{0.1\linewidth}
    \begin{subfigure}{0.4\linewidth}
        \includegraphics[width=\linewidth]{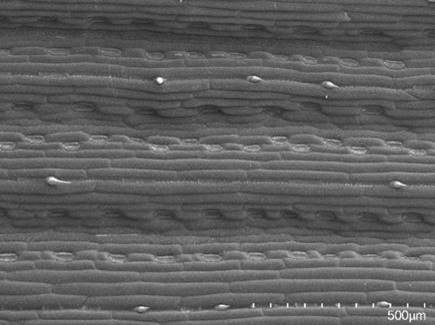}
        \caption{A scanning electron microscope image of winter wheat, showing the type of micro-structures present in our micro-scale dataset.}
        \label{figs:wheatSEM}
    \end{subfigure}
    \caption{Wheat leaves, imaged at two vastly different scales.}
    \label{fig:wheatImages}
\end{figure}

However, compared to the whole leaf, the micro-structures are quite small -- on the order of $10 \mathrm{\mu m}$.
\Cref{figs:wheatSEM} shows an image of the ridges and hairs on a wheat leaf at this scale.
High resolution CT scans are capable of resolving the micro-structures, but the time required for even a small sample, say $2\text{mm}$ of leaf, is many hours, requiring hundreds, if not thousands, of X-ray images to be captured. 
Larger samples require yet longer scanning times to achieve the same resolution, and this quickly reaches a practical limit for an organic sample like a wheat leaf, which over time begins to dry and deform. 
Thus, it is impractical using current scanning technologies to scan a whole leaf with the resolution necessary to capture the micro-structures.

Many approaches have been devised for so-called `multi-scale' or `multi-resolution' surface reconstruction \cite{estellers_robust_2016, hilton_multi-resolution_1997, kuschk_fusion_2013, remondino_3d_2009, siart_fusion_2013, steinbrucker_large-scale_2013, tang_multi-scale_2018, tansky_multi-sensor_2014}.
These approaches usually address the problem of reconstructing a surface from a large dataset using hierarchical representations (multi-scale), or where the density of samples in the dataset is variable (multi-resolution). 
In this work we will also need to deal with large datasets, containing around $100,000$ points, but our problem is more analogous to that of texture synthesis \cite{efros_image_2001}. 
We have access to two methods for scanning plant surfaces at two vastly different scales: we will refer to them as micro- and macro-scale.
The macro-scale data come from an Agri-EPI LemnaTec sensor gantry, and comprise a three-dimensional point cloud of leaf surface points covering an entire plant, from which we have cropped a single leaf. 
The micro-scale data, as already foreshadowed, come from micro-CT scans of only a small patch of leaf surface.  
Thus, our objective is to `paint' the micro-scale reconstruction (the `texture') onto the macro-scale reconstruction, thereby mimicing a data set that would theoretically be obtained by scanning the entire leaf using high-resolution micro-CT.  
Crucially however, our texturing algorithm will generate this representation with only a tiny fraction of the data that would need to be gathered, stored and processed for a full scan. 
As far as we know there are no similar techniques for combining surface reconstructions in this way.

Our method is a new kind of multi-scale model that combines scanned data at vastly different resolutions into one surface reconstruction.
In fact, we will show that the multi-scale reconstruction is a natural extension of an existing method for localising an interpolation to deal with large datasets; the Partition of Unity Method (PUM) \cite{babuska_partition_1997}.
We focus particularly on wheat leaves in this paper.  
Wheat is a vitally important agricultural crop, and it also has particular characteristics that make it both fascinating and challenging to modellers and agricultural scientists alike.
It has a difficult-to-wet leaf due to thin layers of wax ($<200\mathrm{nm}$), but also, as mentioned, ridges and hairs that constrain the shape of droplets when they are deposited and as they evaporate -- \cref{fig:wheatImages} shows these features in detail.
These micro-structures can affect the shape of a sessile droplet, even when compared to a similar difficult-to-wet surface such as Teflon \cite{tredenick_evaporating_2021}.
For these reasons, wheat is a natural target for research on optimising spray technologies \cite{dorr_spray_2016, dorr_towards_2014, dorr_impaction_2015, mayo_simulating_2015, tredenick_evaporating_2021}, which motivates this work.

The rest of this paper is organised as follows: we outline a unifying method for surface reconstruction at both scales in \cref{sec:methods}, show results with a macro-scale model `painted' with micro-scale details in \cref{sec:results}, then discuss the specific contributions and possible extensions of the method in \cref{sec:conclusions}.

\section{Methods}
\label{sec:methods}
\subsection{Datasets}
\begin{figure}
    \centering
    \includegraphics[width=\linewidth]{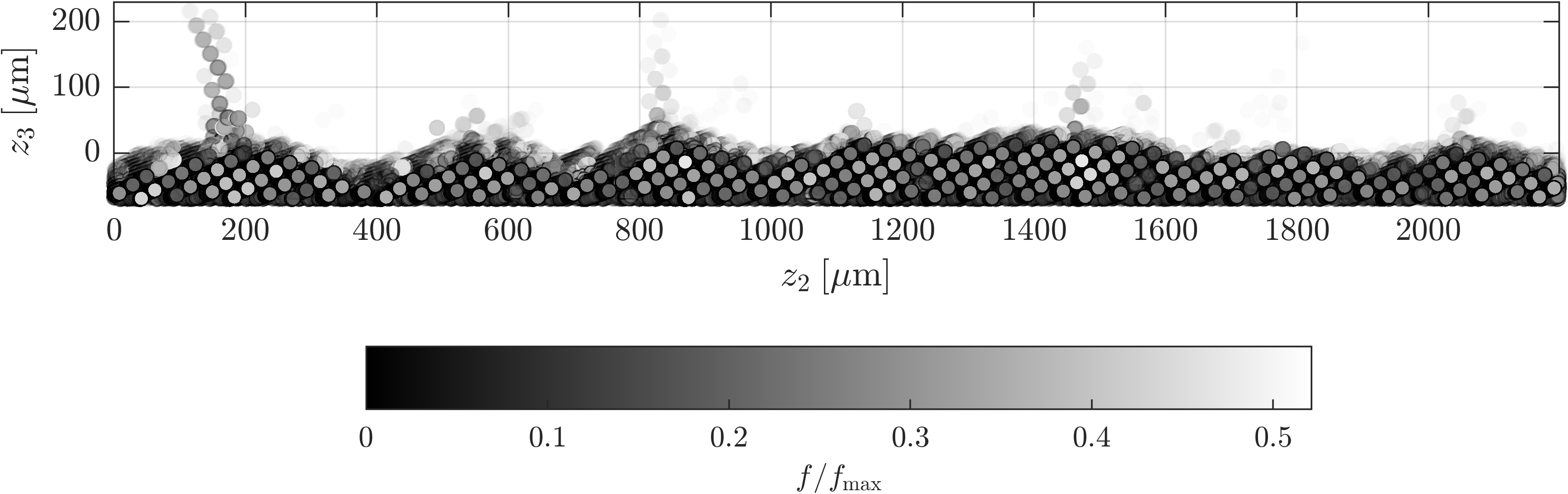}
    \caption{A visualisation of the micro-CT scan dataset. Both the grayscale intensity and the transparency of the points have been mapped to the CT intensity, $f$. The microstructures -- ridges and hairs -- are apparent from this angle. Note that this is a side-on view of a 3D dataset -- see \cref{fig:implicitSurface} for a clearer visualisation of the micro-scale surface.}
    \label{fig:CT_dataset}
\end{figure}

As outlined in \cref{sec:intro}, we will make use of two distinct types of data in this work: micro-CT scan data and point clouds from a range scan.
The micro scale dataset was obtained by scanning a single wheat leaf using high-resolution micro-CT (µCT50, Scanco Medical, Brütisellen, Switzerland). 
The leaf was scanned in air with an isotropic voxel size of 3$\mu \mathrm{m}^3$ at 45 kV, 133 $\mu\mathrm{A}$ and 1.002s sample time. 
892 grayscale  images (slices) were taken, reconstructed from 1500 X-ray projections for each slice. 
The grayscale values correspond to the radio-densities (i.e. X-ray attenuation) of the leaf. 
We consider the radio-densities reconstructed by the CT scanner as indicative of the material present in a voxel: leaf material registers a higher radio-density, while the air registers a lower radio-density.
After appropriately cropping and rotating the dataset, we obtain a set of points (voxel centres) $\vec{z}_i$, together with their radio-density measurements $f_i$.
\cref{fig:CT_dataset} is a visualisation of the micro-CT dataset we use to reconstruct the micro-scale features of the wheat leaf, showing the hairs and periodic ridges that are not apparent at the macro-scale. 

\begin{figure}
    \centering
    \includegraphics[width=0.8\linewidth]{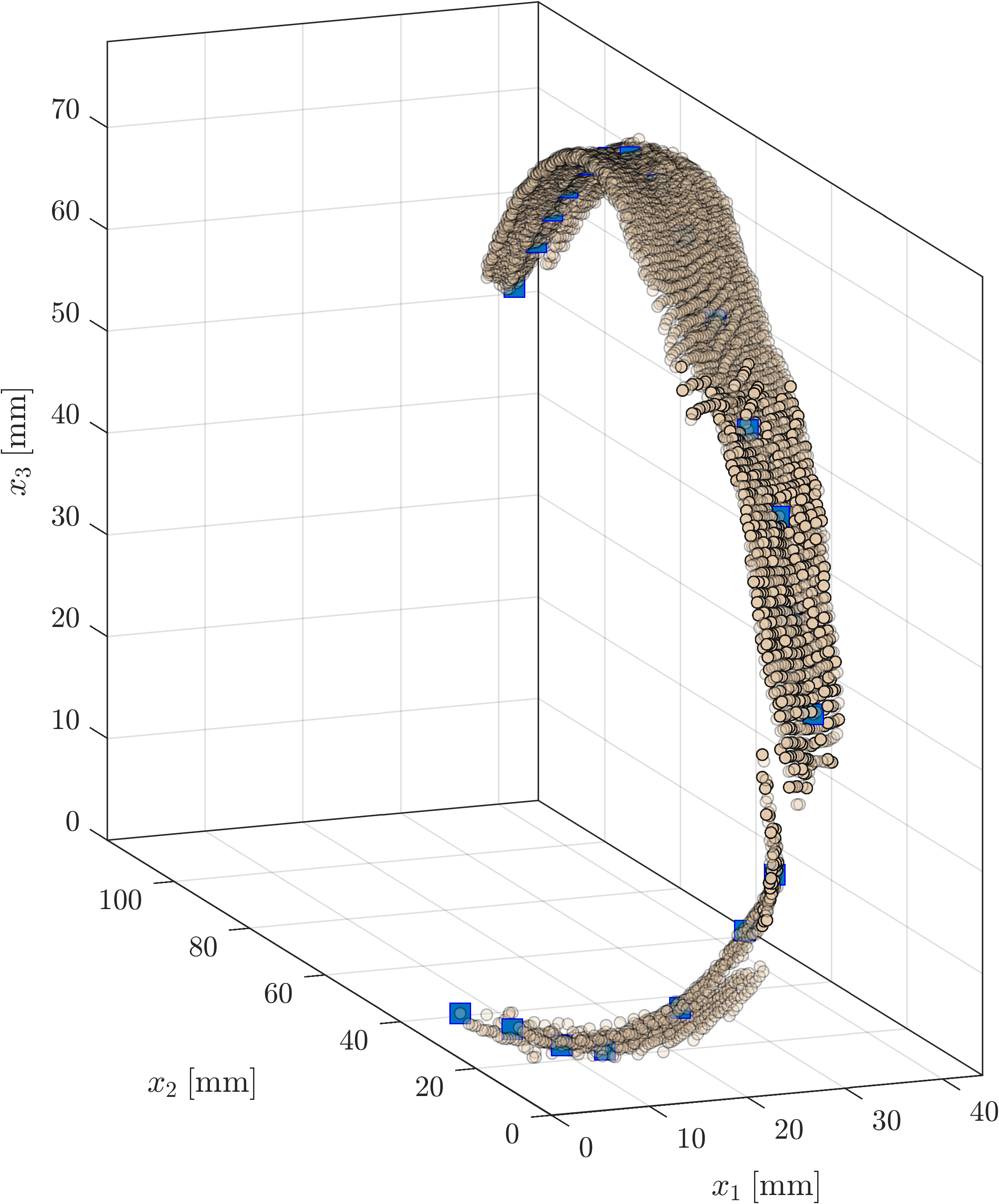}
    \caption{A visualisation of the wheat leaf laser scan dataset. The light brown circles are the scanned points, and the blue squares are the manually selected control points for the medial axis spline. The macroscale features of importance are apparent here: the angle of the surface, the curvature, and the twisting of the leaf.}
    \label{fig:PC_dataset}
\end{figure}

The second type of dataset we use is the laser-scanned point cloud of a whole wheat plant. 
This data was collected by an Agri-EPI LemnaTec sensor gantry laser scanner. 
It comprises simply of a set of points, $\vec{x}_j \in \mathbb{R}^3$, representing one wheat leaf.
\Cref{fig:PC_dataset} shows a visualisation of the datapoints.
We have reconstructed implicit surfaces from precisely this kind of dataset in previous work \cite{whebell_implicit_2021}.

\subsection{Smoothed Radial Basis Function (RBF) interpolation}
\label{sec:methods:RBF}
Our general approach to surface reconstruction is interpolation of scattered data using smoothed radial basis functions (RBF).
We start with a set of points $\{ \vec{x}_j \}_{j=1}^{N} \subset \reals^d$, as well as some measurements at those points, $\{ f_j \}_{j=1}^{N} \subset \reals$.
Interpolants come in two flavours: explicit and implicit.
To fit a height map, the measurements $f_j$ are simply the height of the surface at each sample $\vec{x}_j$. 
Any point on the surface then has the form $[ x_1, x_2, \mathcal{H}(x_1, x_2) ]^T$, where $\mathcal{H}$ is the interpolant.
We use this approach for the macro-scale point cloud data: see \cref{sec:methods:pointCloud}.

For more complex surfaces, not describable by an explicit height map, we may fit an implicit interpolant $\scrF$ such that every point on the surface satisfies $\scrF(\vec{x}) = 0$.
This is the approach we use to fit a surface to the micro-CT data: see \cref{sec:methods:micro}.

Choosing a threshold $f_\text{surface}$ by which to offset the voxel intensities: $f_j \leftarrow f_j - f_\text{surface} \, (j = 1,\dots,N)$, an implicit surface $\scrF(\vec{x}) = 0$ then separates voxel intensities higher than $f_\text{surface}$, where $\scrF(\vec{x}) > 0$, from voxel intensities lower than $f_\text{surface}$, where $\scrF(\vec{x}) < 0$.
It also induces an orientation to the surface: positive on one side, negative on the other.

Whether explicit or implicit, RBF interpolation constructs the function $\scrF(\vec{x})$ as a weighted sum of radial functions $\phi: \reals^+ \rightarrow \reals$, shifted to the datasites $\vec{x}_j$, with the addition of some low-degree polynomial terms, $p_k$:
\[ 
\scrF(\vec{x}) = \sum_{j=1}^N \lambda_j \phi(\| \vec{x} - \vec{x}_j \|) + \sum_{k=1}^n a_k p_k(\vec{x}).
\]
The $\lambda_j$ and $a_k$ are coefficients to be determined.
Imposing interpolation conditions $\scrF(\vec{x}_j) = f_j$ and orthogonality conditions \cite{wahba_spline_1990}:
\[
\sum_{j=1}^N \lambda_j p_k( \vec{x}_j ) = 0
\]
yields an exact RBF interpolant.
The basis functions we use are polyharmonic splines, which are rotationally invariant functions minimising a particular measure of the interpolant's curvature: the thin-plate penalty functional of order $m$ on $\mathbb{R}^d$ \cite{wahba_spline_1990}:
\begin{equation} \label{energy}
	J_m^d(\mathcal{F}) = 
	\sum_{\alpha_1+\dots+\alpha_s=m} \:
	\frac{m!}{\alpha_1! \dots \alpha_d!}
	\int_{\mathbb{R}^d}
	\left(
	\frac{ \partial^m \mathcal{F} }{ \partial x_1^{\alpha_1}\dots\partial x_d^{\alpha_d} }
	\right)^2
	\: \mathrm{d}\vec{x},
\end{equation}
which is the sum of square-integrated $m$-th order mixed partial derivatives of $\mathcal{F}$.
The general form of the polyharmonic spline of order $m$ in dimension $d$ is \cite{duchon_splines_1977, wahba_spline_1990}:
\begin{align*}
	E(\tau) &= 
	\begin{cases}
		\theta_m^d |\tau|^{2m-d} \log|\tau|, & \text{$2m-d$ even}, \\
		\theta_m^d |\tau|^{2m-d}, & \text{otherwise},
	\end{cases}
	\\
	\text{with} \quad \theta_m^d &= 
	\begin{cases}
	\frac{(-1)^{d/2+1+m}}{ 2^{2m-1} \pi^{d/2} (m-1)! (m-d/2)! },
	& \text{$2m-d$ even}, \\
	\frac{\Gamma(d/2-m)}{ 2^{2m} \pi^{d/2} (m-1)! }, 
	& \text{otherwise}.
	\end{cases}
\end{align*}
In this work we use the popular thin plate spline $\phi(r) = r^2 \log r$ in 2D, and $\phi(r) = r^3$ in 3D.

Not only do these splines produce interpolants of minimal curvature, they can also be used to fit approximate, smoothed interpolants to noisy data.
We can fit a smoothed surface to noisy data by solving the penalised least-squares problem \cite{wahba_spline_1990}:
\[
\argmin_{\vec{\lambda}, \vec{a}} 
\frac{1}{N} \sum_j \left[ \scrF(\vec{x}_j) - f_j \right]^2 + \rho J_m^d(\scrF),
\]
where $J$ is the same thin-plate penalty functional as in equation \cref{energy}.
This minimisation problem can be cast as the linear system \cite{wahba_spline_1990}:
\begin{equation}
    \left( A + \rho N / \theta^d_m I \right) \vec{\lambda} + P \vec{a} = \vec{f}, 
    \label{eq:linSys}
\end{equation}
where:
\begin{gather*}
    A_{ij} = \phi( \| \vec{x}_i - \vec{x}_j \|_2 ), \\
    P_{ik} = p_k(\vec{x}_i), \\
    \vec{\lambda} = [\lambda_1, \dots, \lambda_N]^T, \quad
    \vec{a} = [a_1, \dots, a_n]^T, \\
    i,j = 1, \dots, N, \quad
    k = 1, \dots, n,
\end{gather*}
subject to:
\begin{gather*}
    P^T \vec{\lambda} = \vec{0},
\end{gather*}
where $I$ is the $N\times N$ identity matrix.
The constant $\theta^d_m$ is a normalisation factor dependent on the particular choice of polyharmonic spline.
For our splines this constant takes the values $8\pi$ in 2D, and $96\pi$ in 3D.
The smoothing parameter $\rho$ controls the smoothness of the inexact interpolant, with $\rho = 0$ corresponding to an exact fit. 
We choose $\rho$ qualitatively from experience with the datasets, and we find that the results are not particularly sensitive to its value when varied within an order of magnitude.
More detail on polyharmonic splines in general is given in \cite{wahba_spline_1990}, and specifically in the context of implicit surface reconstruction in \cite{whebell_implicit_2021}.

\subsection{Partition of Unity Method (PUM)}
\label{sec:methods:PUM}
The coefficient matrix in \cref{eq:linSys} is dense. 
For a dataset with $N$ points, and an interpolant with $n$ polynomial terms, this is an $(N+n)\times(N+n)$ linear system, which is infeasible for large datasets such as ours ($N \approx 100,000$ points).
Using the partition of unity method, we can introduce locality to the interpolation by partitioning the domain into overlapping subdomains, which overcomes this limitation \cite{franke_locally_1977, wendland_computational_2006}.
We use hyperspheres $\{\Omega_i\}_{i=1}^m$ defined by a centre $\vec{c}_i$ and a radius $r_i$.
We then independently fit local interpolants, $\scrF_i$, to the data contained in each subdomain.
The global interpolant is a sum of the local interpolants, weighted by distance to subdomain centres:
\begin{gather*}
    \scrF(\vec{x}) = \sum_{i=1}^M w_i(\vec{x}) \scrF_i(\vec{x}).
\end{gather*}
The weight function we use is a compactly supported Wendland function \cite{wendland_piecewise_1995}:
\begin{equation*}
    w(\tau) = 
    \begin{cases}
        (1 - \tau)^4 (4\tau + 1), & 0 \leq \tau \leq 1, \\
        0 & 1 < \tau,
    \end{cases}
\end{equation*}
which is maximised at a subdomain's centre, $\vec{c}_i$, and vanishes outside the subdomain.
The weights at any one point are normalised to sum to unity using Sheppard's method \cite{shepard_two-dimensional_1968}:
\begin{gather*}
    w_i(\vec{x}) = \frac{ 
        w( \| \vec{x} - \vec{c}_i \| / r_i )
    }{
        \sum_k w( \| \vec{x} - \vec{c}_k \| / r_k )
    }.
\end{gather*}

If subdomains are chosen such as to keep the number of datapoints in each subdomain ($N/M$) constant, then the complexity of fitting the local interpolants is $\mathcal{O}( M (N/M)^3 ) = \mathcal{O}( N )$.
We can satisfy this condition by adaptively choosing the subdomains according to the density of the datapoints using, for example, the octree subdivision method given in \cite{whebell_implicit_2021}.
In \cref{sec:results} we will show that the PUM not only reduces the complexity of the interpolation problem, it can also `blend' together different, similar, surfaces at their edges.

\subsection{Surface reconstruction from point cloud (macro-scale)}
\label{sec:methods:pointCloud}
\begin{figure}
    \centering
    \includegraphics[width=0.9\textwidth]{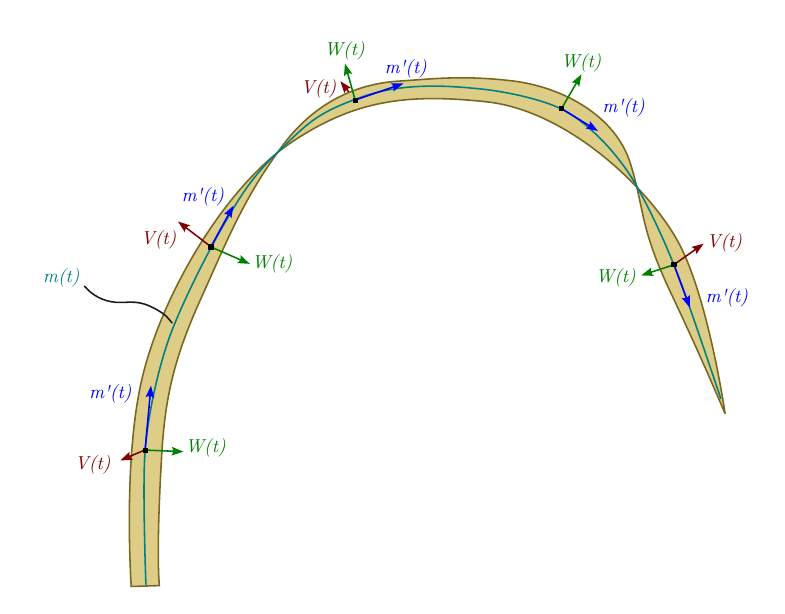}
    \caption{A schematic of the locally orthogonal coordinate system following the leaf surface. The medial axis of the leaf is modelled as a spline $\vec{m}(t)$, the tangent of which is the first coordinate, $\vec{m}'(t)$. The next coordinate is the approximate surface normal $\vec{V}(t)$ and the final coordinate is $\vec{W}(t)$, the cross product of the first two.}
    \label{fig:orthoCoords}
\end{figure}

To build an interpolation for the macro-scale point cloud dataset shown in \cref{fig:PC_dataset}, we observe that the wheat leaves we seek reconstructions of are reasonably flat, locally.  
Hence, an explicit interpolation approach is feasible.  
First, we manually choose some control points from the point cloud that lie along the medial axis of the leaf blade.  
Fitting a spline through these control points defines a parmeterisation following the curvature of the leaf, from which we may build a full local coordinate system for the leaf surface, which we refer to as the leaf coordinate system: see \cref{sec:orthoCoords} for details on the construction. 
A mathematical transform to and from the world coordinates, $\vec{y} = \mathcal{T}(\vec{x})$ then follows, in a similar manner to the approach used in \cite{kempthorne_surface_2014}.  
This coordinate system is illustrated in \cref{fig:orthoCoords}.  
Hence by transforming the points $\vec{x}_j$ to leaf coordinates $\vec{y}_j$, we may reconstruct the macro-scale leaf surface with a smoothed, explicit RBF interpolant $\mathcal{H}$ as detailed in \cref{sec:methods:RBF}:
\[ \mathcal{H}(y_1, y_2) \approx y_3. \]

\subsection{Surface reconstruction from CT scan (micro-scale)}
\label{sec:methods:micro}
\begin{figure}
    \centering
    \includegraphics[width=0.7\linewidth]{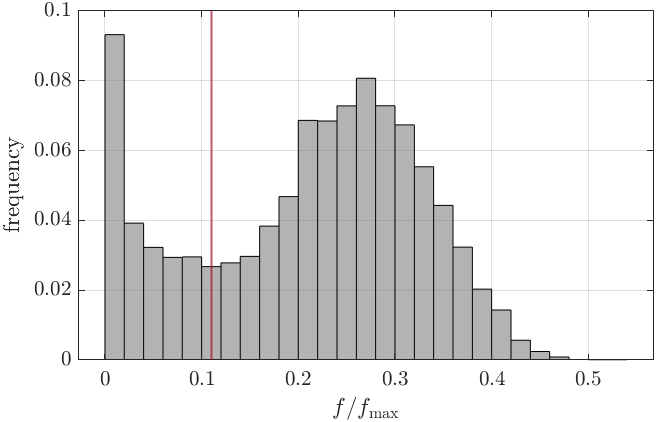}
    \caption{A histogram of relative voxel intensities in the micro-CT dataset (zero-valued voxels excluded). The red line shows a good choice for the threshold $f_\text{surface}$ to separate leaf material from the air.}
    \label{fig:CT_histogram}
\end{figure}

When it comes to the micro-scale CT data, the hairs on the wheat leaf are an essential feature to be captured in the representation (see, for example, \cref{fig:implicitSurface}).  
This necessitates an implicit interpolant, since there is no simple coordinate transform that would represent these features as an explicit height map.  
We fit an implicit, smoothed RBFPUM interpolant $\mathcal{F}(\vec{z}) = f_\text{surface}$ to the datapoints $(\vec{z}_i, f_i)$, as detailed in \cref{sec:methods:RBF,sec:methods:PUM}.  
The threshold value $f_\text{surface}$ is determined by generating a histogram of voxel intensities in the CT scan dataset see \cref{fig:CT_histogram}.
Simply choosing an intensity between the two peaks, corresponding to air and leaf material, provides a suitable threshold value; for example the relative value $0.11$ for the dataset in \cref{fig:CT_histogram}.
For the wheat leaf, we found that higher values in this range can leave voids inside the leaf where there may have been air or water pockets, so for this work we have used $f_\text{surface} = 0.03$.

\subsection{Multi-scale reconstruction}
\label{sec:methods:multiscale}
To combine the macro- and micro-scale features into a single multi-scale reconstruction, we start with the leaf coordinate system for the macro-scale point cloud. 
We will, in essence, tile the plane that the leaf lies on in the leaf coordinate system with copies of our micro-scale surface patch. 
To achieve this outcome, we consider the micro-scale interpolant $\scrF(\vec{z})$ and its PUM subdomains $\Omega_i$. 
We take care to appropriately rotate the micro-scale data before fitting to ensure that the leaf aligns with the $(z_1, z_2)$ plane, to simplify the next step. 
The subdomain $\Omega_i$ and its associated local interpolant $\scrF_i$ are then copied and displaced along the $(z_1, z_2)$ plane by exactly the width of the micro-scale segment. 
We can repeat this process, copying and displacing subdomains, until we have made enough copies to cover the flat leaf point cloud in the $(y_1, y_2)$ plane, denoting the number of copies $N_q$. 
Where one copy of the leaf segment meets another, the PUM will continuously blend the local interpolants together despite the scanned surface not being truly periodic.
Note, though, that in the implementation, we need not actually copy the datapoints or the local interpolant in memory. 
We simply record the displacement vector, which for the $q$-th copy of the $i$-th subdomain we denote $\vec{s}_i^q$. 

Finally, the multi-scale interpolant at point $\vec{x}$ can be defined in world coordinates as:
\[
    \mathcal{G}(\vec{x}) = 
    \sum_{i=1}^{m} 
    \sum_{q=1}^{N_q}
    w_i^q(\vec{x}) \scrF_i(\mathcal{T}^*(\vec{x}) - \vec{s}_i^q),
\]
where
\[
    w_i^q(\vec{x}) = 
    \frac{ 
        w \left( \left\| \mathcal{T}^*(\vec{x}) - \vec{s}_i^q - \vec{c}_i \right\| / r_i \right)
    }{
        \sum_{k=1}^m \sum_{q=1}^{N_q} w \left( \left\| \mathcal{T}^*(\vec{x}) - \vec{s}_k^q - \vec{c}_k \right\| / r_k \right)
    }
\]
is the new normalised weight function to complete the partition of unity.
The function $\mathcal{T}^*(\vec{x})$ involves transforming a point from world coordinates to leaf coordinates, and also incorporates the macro-scale height map. 
Denote $\mathcal{T}(\vec{x}) = \vec{y} = [y_1, y_2, y_3]^T$. 
Then:
\[
\mathcal{T}^*(\vec{x}) = [ y_1, y_2, y_3 - \mathcal{H}(y_1, y_2) ]^T.
\]

This procedure is summarised in \cref{alg:evalMultiScale}.

\begin{algorithm}
    \caption{Evaluating the multi-scale interpolant}
    \label{alg:evalMultiScale}
    \begin{algorithmic}
        \Require query point $\vec{x}$, weight function $w$, local interpolants $\{\mathcal{F}_i\}_{i=1}^m$, world-to-leaf transform $\mathcal{T}$, macro-scale interpolant $\mathcal{H}$, displacement vectors $\{\vec{s}_i^q\}_{i=1, q=1}^{i=m, q=N_q}$, subdomain centres $\{\vec{c}_i\}_{i=1}^m$, subdomain radii $\{r_i\}_{i=1}^m$
        \Ensure $\mathcal{G}(\vec{x})$
        \Function{$\mathcal{T}^*$}{$\vec{x}$}
            \State $[y_1, y_2, y_3]^T \gets \mathcal{T}(\vec{x})$
            \State $\vec{y}^* \gets [ y_1, y_2, y_3 - \mathcal{H}(y_1, y_2) ]^T$
            \State \Return $\vec{y}^*$
        \EndFunction
        \State $g \leftarrow 0$ 
        \Comment{Function value}
        \State $w_\text{sum} \leftarrow 0$
        \Comment{Sheppard's method sum}
        \State $\vec{z} \gets \mathcal{T}^*(\vec{x})$
        \For{$i=1$ to $m$} \Comment{for each subdomain}
            \For{$q=1$ to $N_q$} \Comment{for each copy}
                \State $\tau \leftarrow \| (\vec{z} - \vec{s}_i^q - \vec{c}_i) / r_i \|_2$
                \Comment{normalised distance to copy $q$ of subdomain $i$}
                \State $w_\text{sum} \leftarrow w_\text{sum} + w(\tau)$
                \State $g \leftarrow g + w(\tau) \mathcal{F}_i(\vec{z})$
                \Comment{only compute this if $\tau < 1$}
            \EndFor
        \EndFor
        \State $g \leftarrow g / w_\text{sum}$
        \State \Return g
    \end{algorithmic}
\end{algorithm}

\section{Multi-scale wheat leaf surface reconstruction}
\label{sec:results}
\subsection{Computational details}
We have implemented the method outlined in \cref{sec:methods} in the programming language \texttt{julia} \cite{bezanson_julia_2017}, using several open-source packages (including, but not limited to \cite{barbary_dierckxjl_2022,seth_bromberger_juliagraphslightgraphsjl_2017, k_mogensen_optim_2018,fairbanks_james_juliagraphsgraphsjl_2021}), and made the code \href{https://github.com/rwhebell/Whebell2023_MultiscaleWheat}{available online}.

Using the partition of unity method, the computational complexity of fitting the interpolants scales linearly ($\mathcal{O}(N)$) with the number of datapoints, and the marginal cost of evaluating the interpolant at one point is $\mathcal{O}(1)$.
In \cref{tab:computationTimings} we report the relative timings of different stages of the method, measured on a machine with an AMD Ryzen 2700 (8 cores, 3.2GHz).

\begin{table}
    \renewcommand{\arraystretch}{1.4}
    \centering
    \caption{Relative timings of parts of the method as implemented in \texttt{julia}.}
    \begin{tabular}{p{0.4\linewidth}cc}
        \toprule
        Operation & Scale & Time \\
        \midrule
        Fit micro-scale interpolant & 98k points, 128 subdomains & 16s \\
        Fit leaf coordinate system & 19 points & 6ms \\
        Use world to leaf transform & 6.3k points & 1.0s \\
        Fit macro-scale interpolant & 6.3k points, 4 subdomains & 5.5s \\
        Using $\alpha$-shape to exclude exterior points from evaluation & 25M points & 15 mins \\
        Evaluating implicit function at interior points & 1.9M points & 3 mins \\
        Polygonizing surface with marching tetrahedra & 130M tetrahedra & 6.5 mins \\
        Sample height map & 8.9M points & 12 mins \\
        Transform triangles back to world coords & 8.9M points & 71s \\
        \bottomrule
    \end{tabular}
    \label{tab:computationTimings}
\end{table}

\subsection{Locally orthogonal coordinate system}
\begin{figure}
    \centering
    \begin{subfigure}{0.45\linewidth}
        \centering
        \includegraphics[width=0.7\linewidth]{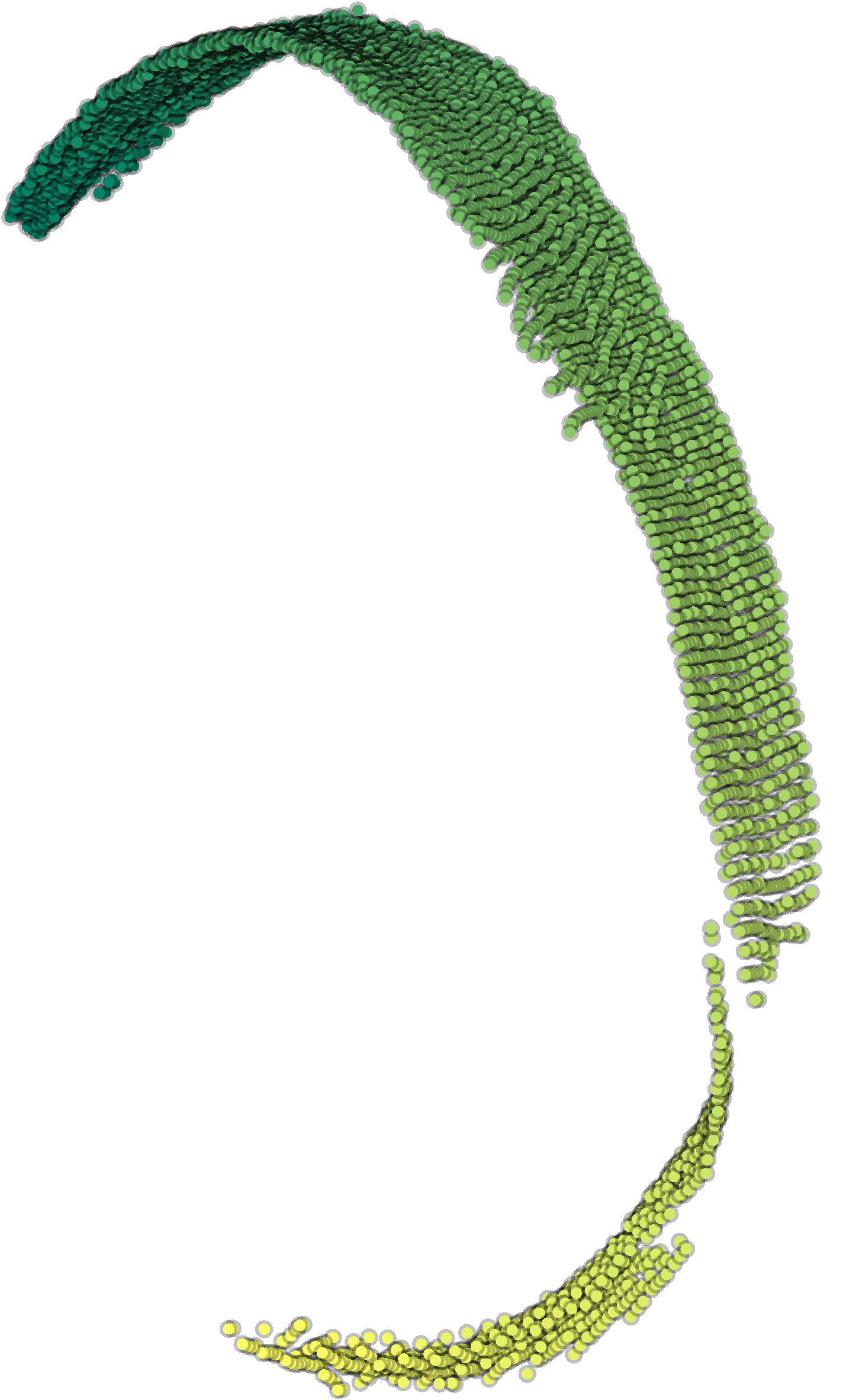}
        \caption{The point cloud dataset from the laser scan of a wheat plant, in world coordinates. The points are coloured by their first coordinate in the orthogonal coordinate system, representing how far along the medial axis a point lies.}
        \label{fig:curlyPointCloud}
    \end{subfigure}
    \hspace{0.08\linewidth}
    \begin{subfigure}{0.45\linewidth}
        \centering
        \includegraphics[width=\linewidth]{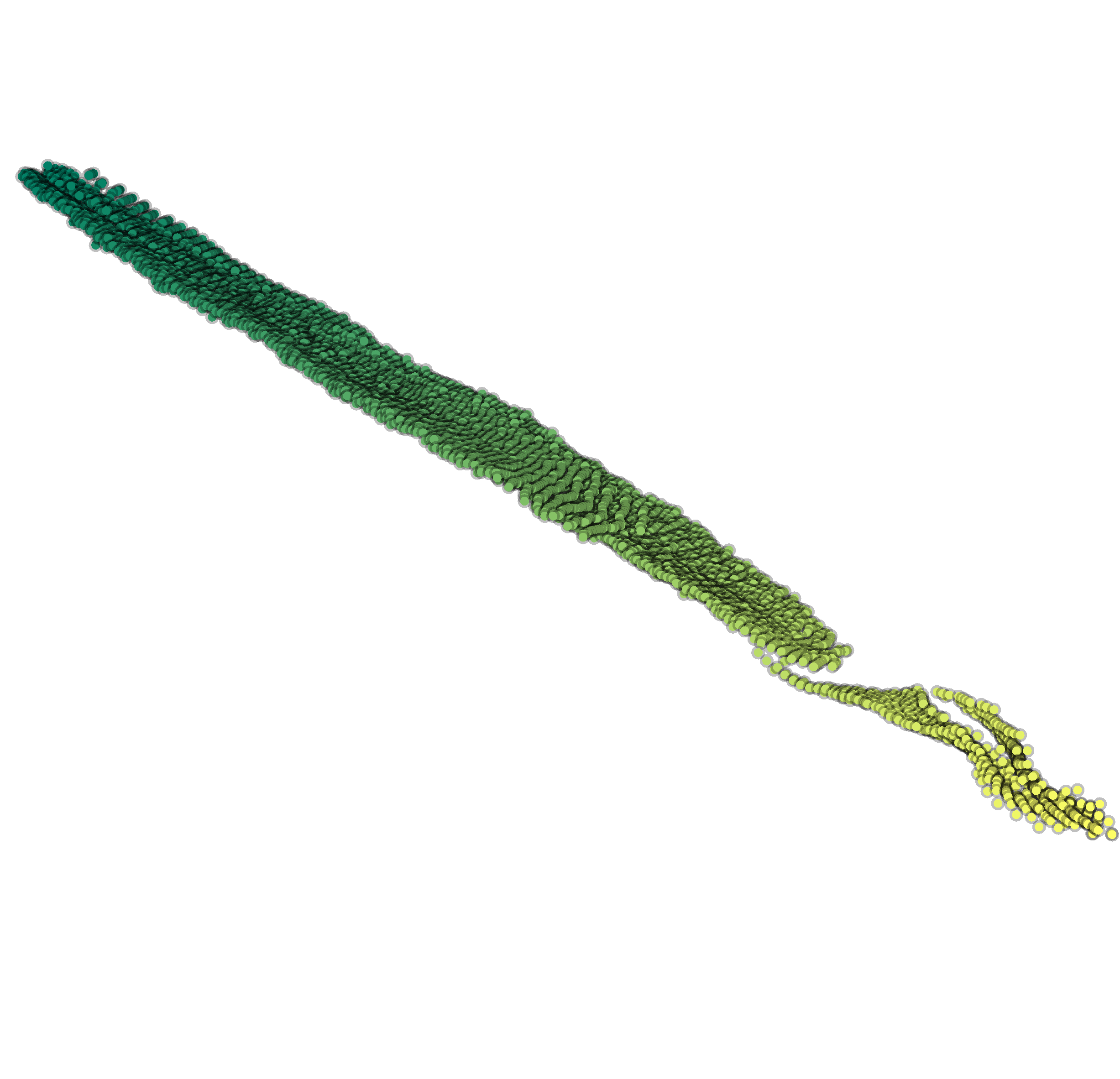}
        \caption{The same point cloud as in \cref{fig:curlyPointCloud}, transformed with $\mathcal{T}$ to leaf coordinates.}
        \label{fig:flatPointCloud}
    \end{subfigure}
    \caption{The effect of the coordinate transform $\mathcal{T}$.}
    \label{fig:world2leaf}
\end{figure}

\begin{figure}
    \centering
    \begin{subfigure}{0.45\linewidth}
        \centering
        \includegraphics[width=\linewidth]{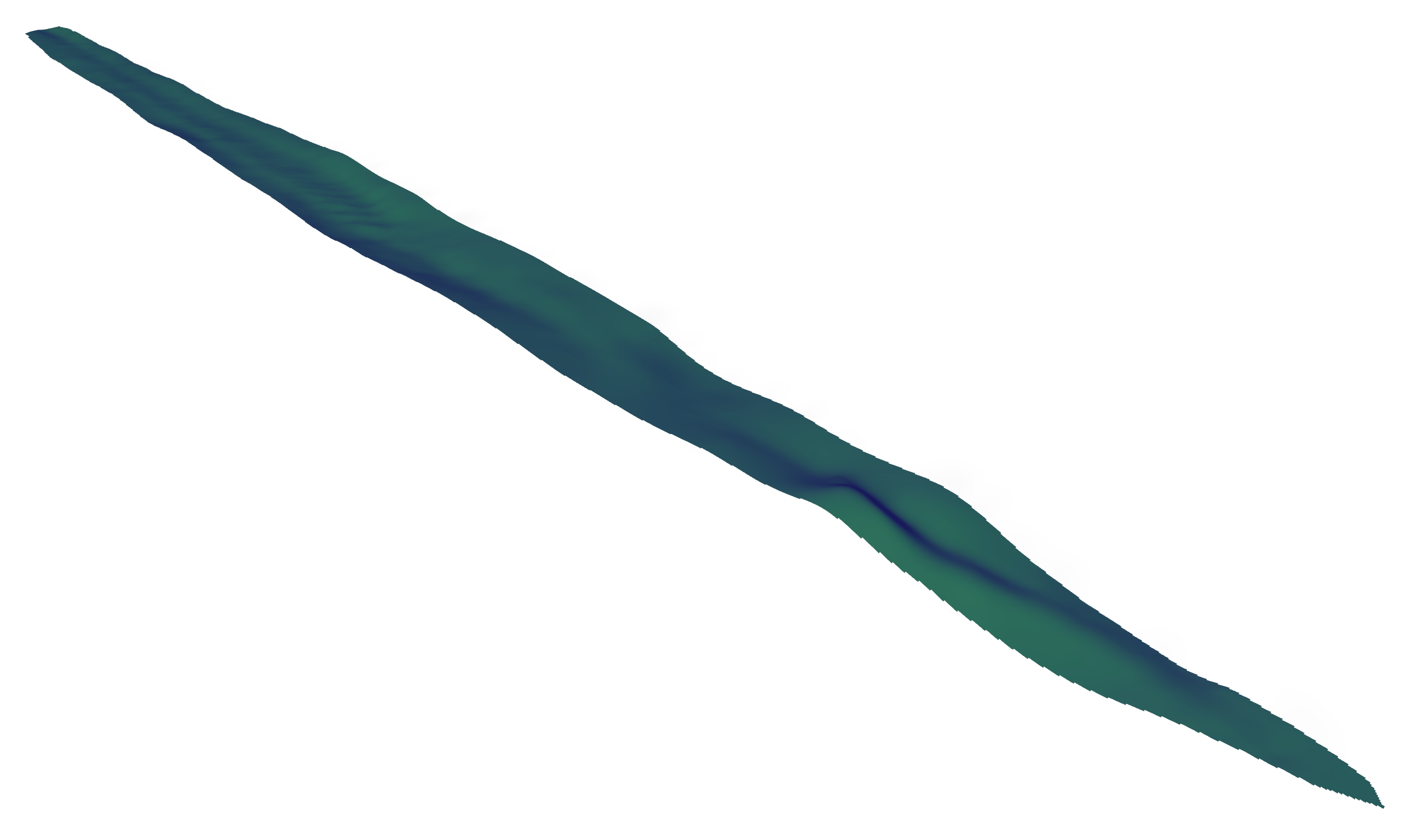}
        \caption{An explicit surface reconstruction from the point cloud, shown in leaf coordinates.}
        \label{fig:flatLeafHeightMap}
    \end{subfigure}
    \hspace{0.05\linewidth}
    \begin{subfigure}{0.45\linewidth}
        \centering
        \includegraphics[width=0.7\linewidth]{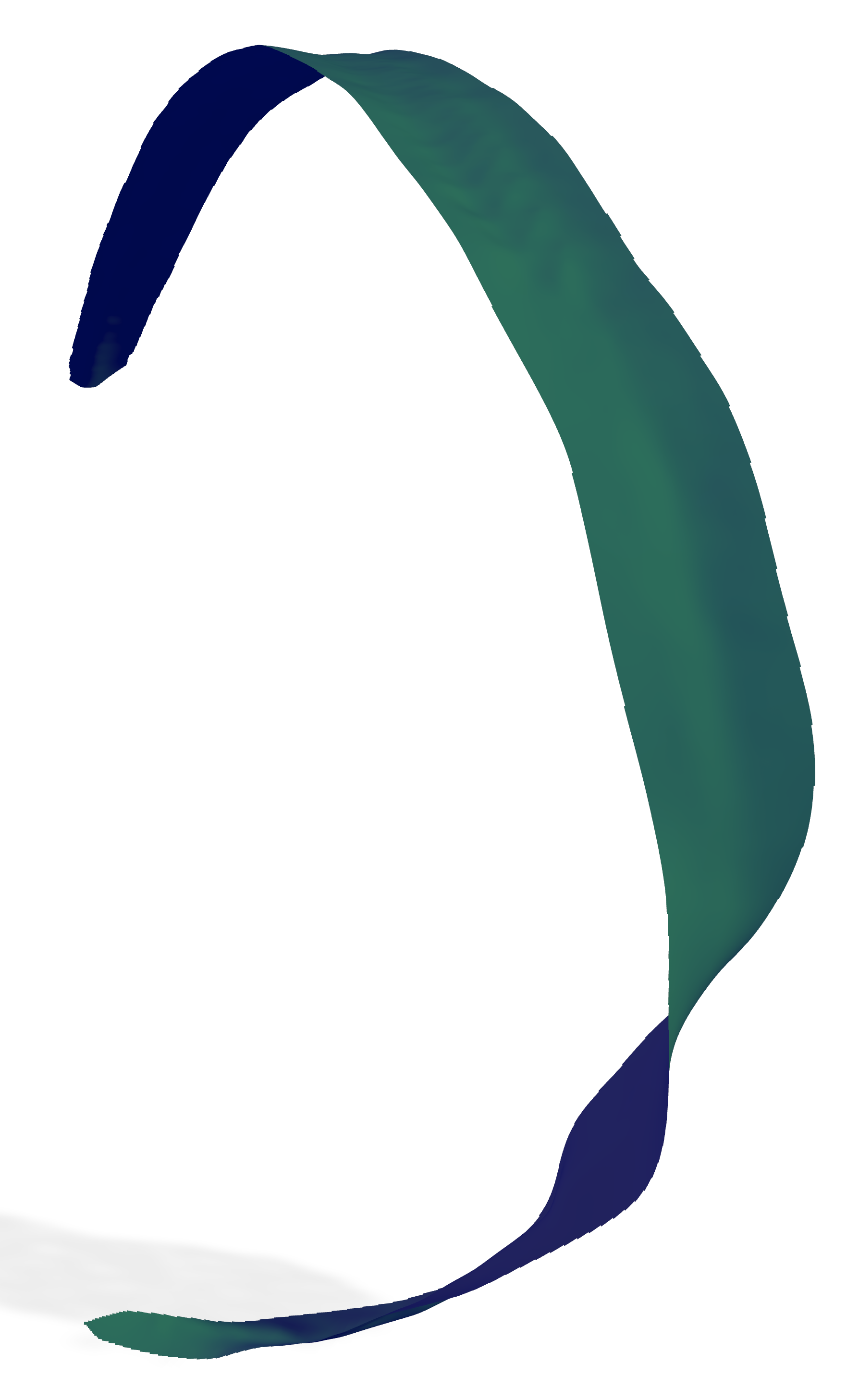}
        \caption{An explicit surface reconstruction from the point cloud, shown in world coordinates.}
        \label{fig:curlyLeafHeightMap}
    \end{subfigure}
    \caption{Demonstrating the explicit surface reconstruction for the macro-scale dataset.}
\end{figure}

\cref{fig:world2leaf} illustrates the effect of the coordinate transform $\mathcal{T}$ on the macro-scale point cloud dataset. 
In the new locally orthogonal coordinate system following the leaf, the surface is single-valued, permitting an explicit interpolant to be fit to the transformed data. 
This surface was fit to around 6,300 points, requiring 5.5 seconds of compute time.
This explicit interpolant is shown in the leaf coordinates in \cref{fig:flatLeafHeightMap} and in the world coordinates in \cref{fig:curlyLeafHeightMap}.
We observe that the simple and efficient parameterisation of the surface by means of the locally orthogonal coordinate system leads to a simple, explicit, surface reconstruction, similar to the approach previously taken in \cite{kempthorne_surface_2014}.
Next we will demonstrate how this macro-scale surface reconstruction not only captures the larger features of the leaf (shape, size, orientation), but also allows us to tile the leaf with copies of the micro-scale surface model. 

\subsection{Partition of unity blending edges}
\begin{figure}
    \centering
    \begin{subfigure}{0.45\linewidth}
        \centering
        \includegraphics[width=\linewidth]{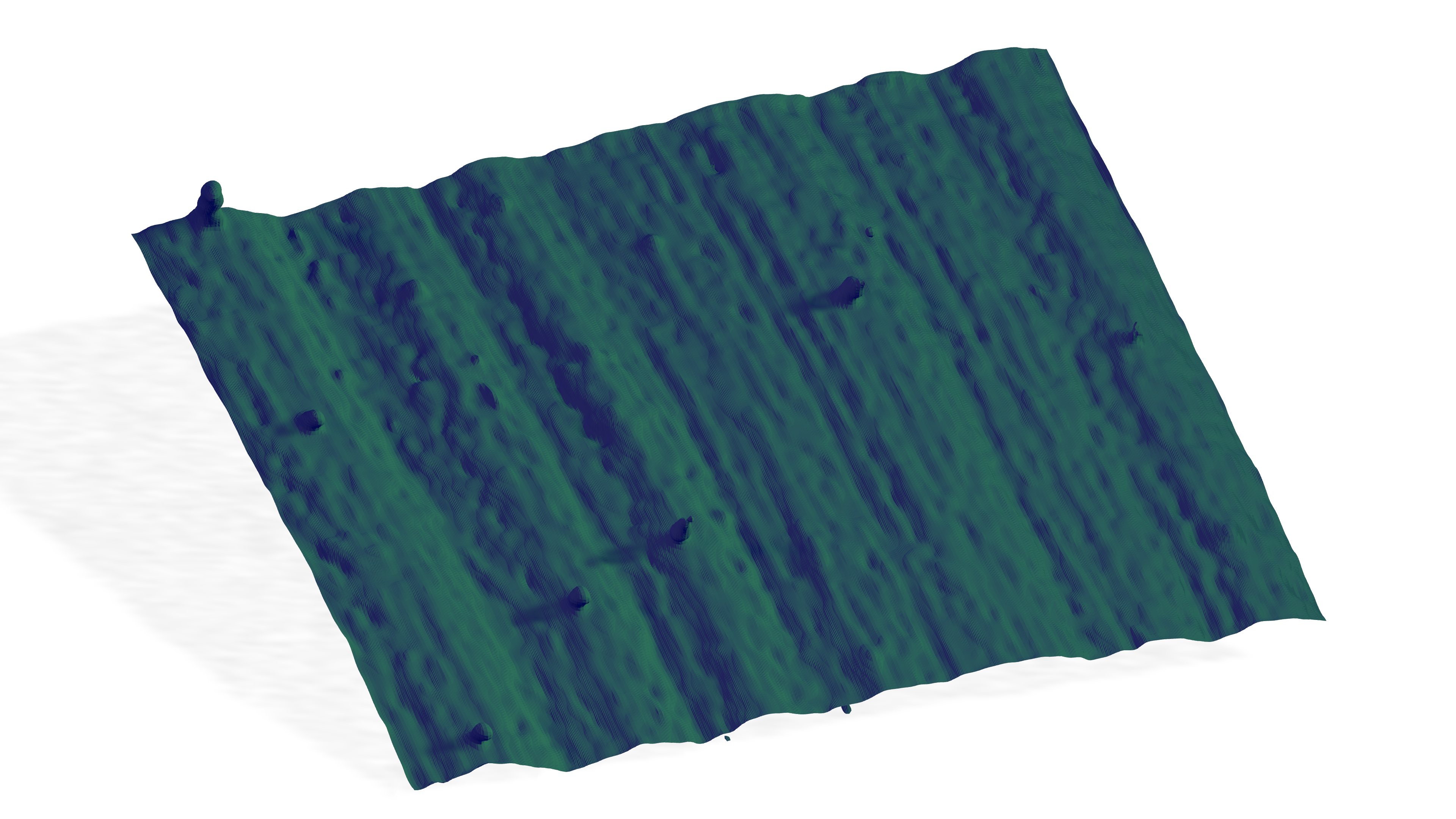}
        \caption{One section of the multi-scale wheat leaf reconstruction, showing just one copy of the micro-scale implicit surface, approximately 2mm wide. Some hairs are visible on the surface (e.g., top left).}
        \label{fig:implicitSurface}
    \end{subfigure}
    \hspace{0.05\linewidth}
    \begin{subfigure}{0.45\linewidth}
        \centering
        \includegraphics[width=\linewidth]{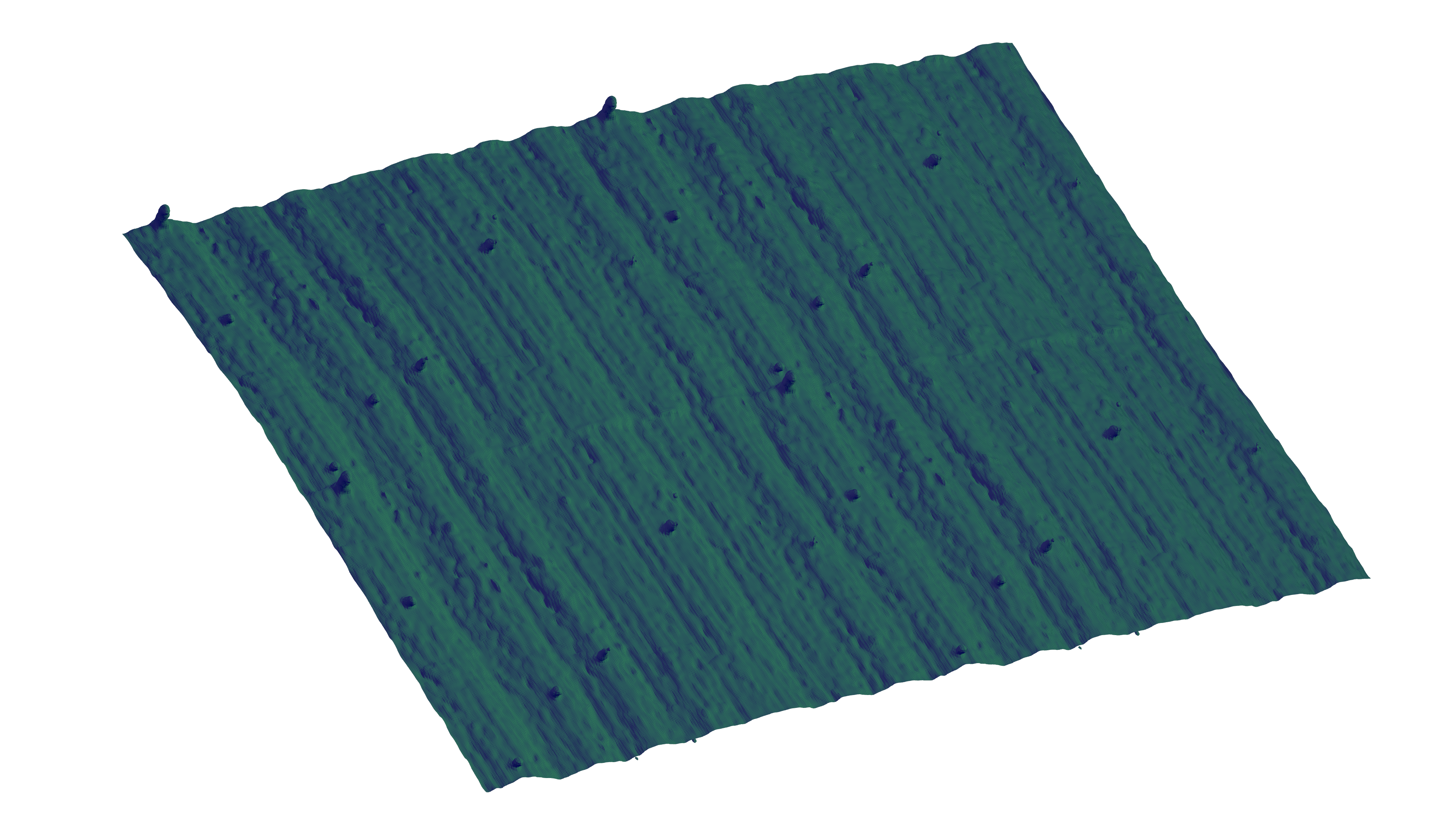}
        \caption{A small section of the multi-scale surface reconstruction, showing four copies of the micro-scale surface, shifted and blended together at the edges.}
        \label{fig:PUMblending}
    \end{subfigure}
    \caption{Demonstrating the Partition of Unity Method `blending' the surface together at edges where distinct copies of the micro-scale interpolant meet.}
\end{figure}

Here we show the partition of unity method blending together edges of the small patch of micro-scale surface, to tile the virtual leaf surface. 
The micro-scale surface reconstruction is fit to around 98,000 points in 16 seconds as described in \cref{sec:methods:micro}.
\cref{fig:implicitSurface} shows the implicit surface for reference, and \cref{fig:PUMblending} shows an implicit surface constructed by shifting three copies of each of the partition of unity subdomains as described in \cref{sec:methods:PUM}.
A slight `seam' is visible through the middle of this surface (running left-to-right in the figure) as a result of the surface not being truly periodic.
This artefact persists in the multi-scale reconstruction, but the surface nonetheless remains continuous, and the seam is small compared to the other micro-scale model's features.  Hence it is not of concern when it comes to running droplet simulations on such a surface (as is the motivation for the present work).
This natural tiling of micro-scale features achieved through blending partition of unity interpolants is key to the multiscale reconstruction to be described in the next section.

\subsection{Multi-scale wheat leaf reconstruction}
With the ability to blend copies of the micro-scale surface together at their edges, we can now readily tile the plane in flattened leaf coordinate space with translated copies of the single patch.
Conceptually, we arrange the copies in a simple grid. 
Computationally, there is no actual duplication of any data; only the offset values of each tile need to be stored, along with just one copy of the actual surface data.
To maintain the shape of the outline of the leaf, we use an $\alpha$-shape \cite{edelsbrunner_three-dimensional_1992} of the point cloud in flat leaf coordinates to define the domain for the implicit interpolant.  
As a generalisation of the convex hull to non-convex volumes, the $\alpha$-shape is realised by building a Delaunay triangulation and removing the largest elements (those with the largest circumradii, up to a user-defined tolerance). 
When realising the implicit surface using marching cubes, triangles with vertices outside the $\alpha$-shape are discarded.  
This method, and a more thorough description of $\alpha$-shapes, can be found in \cite{whebell_implicit_2021}.

\begin{figure}
    \centering
    \includegraphics[width=0.8\linewidth]{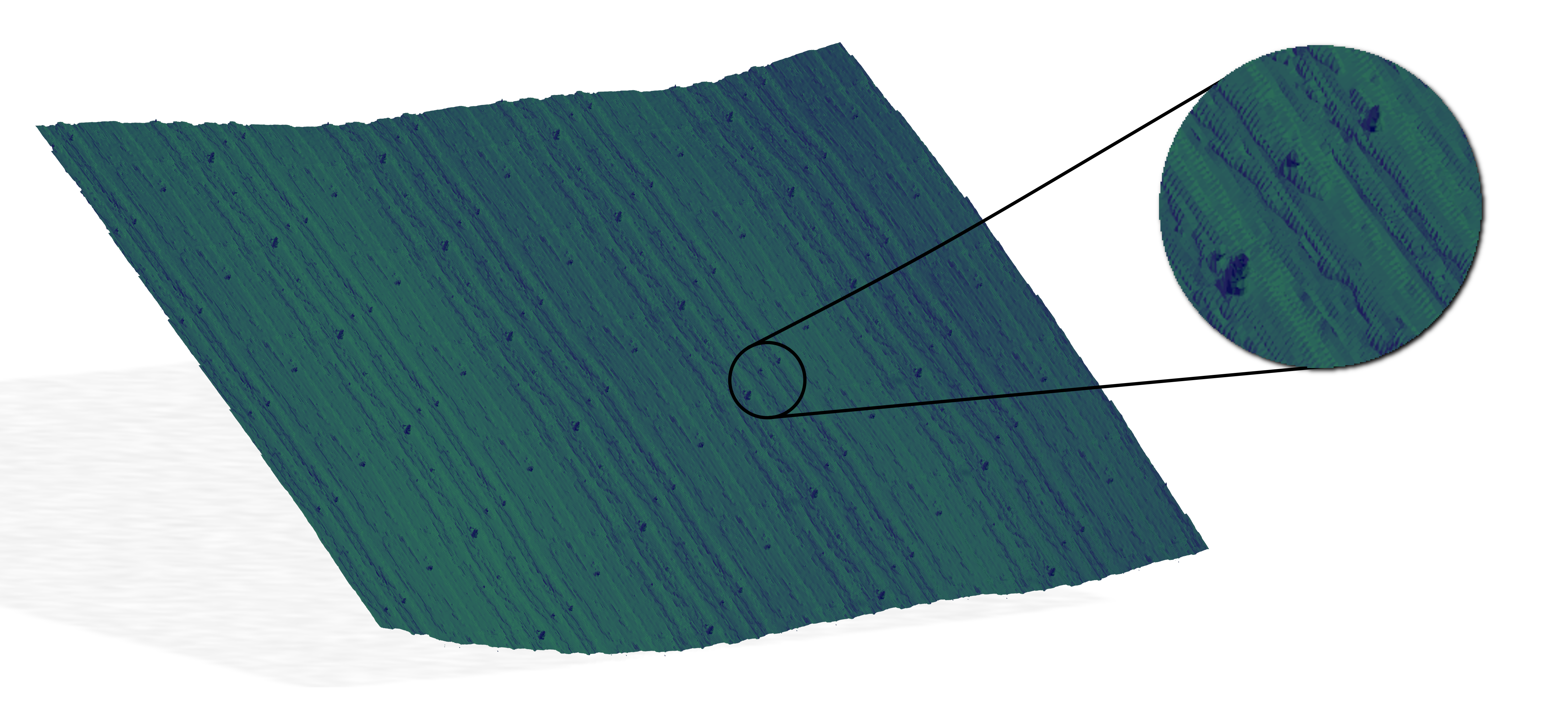}
    \caption{A cross-section of the multi-scale wheat leaf surface reconstruction, approximately 1cm wide. Inset: a close-up of some micro-scale hairs featured on the surface.}
    \label{fig:wheatCrossSection}
\end{figure}

\begin{figure}
    \centering
    \includegraphics[width=\linewidth]{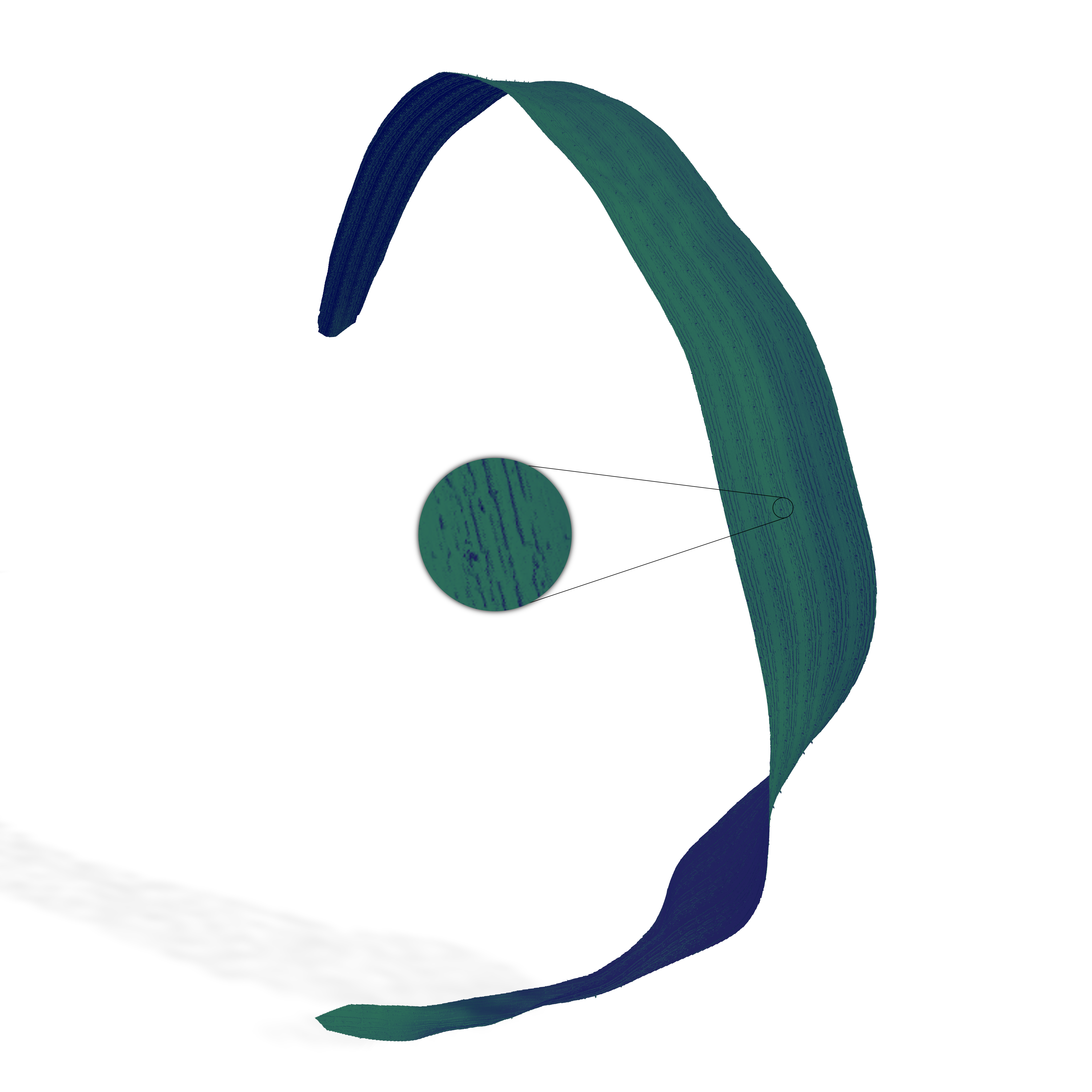}
    \caption{A visualisation of the whole multi-scale wheat leaf surface reconstruction. Inset: a close-up showing some micro-scale details still present in the multi-scale reconstruction.}
    \label{fig:wholeWheatLeaf}
\end{figure}

\cref{fig:wheatCrossSection} shows a cross-section of the wheat leaf reconstruction, showcasing both macro-scale curvature from the height map and micro-scale detail from the implicit surface. 
\cref{fig:wholeWheatLeaf} shows a visualisation of the whole reconstruction, in the shape of a wheat leaf with micro-scale textures on the surface.
Every figure in this and the last subsection (\cref{fig:implicitSurface}, \cref{fig:PUMblending}, \cref{fig:wheatCrossSection}, and \cref{fig:wholeWheatLeaf}) was generated using the same surface reconstruction, evaluated at different resolutions over different subsets of the domain.

\section{Conclusions}
\label{sec:conclusions}
In this paper, we have described and demonstrated a unifying method for reconstructing a macro-scale surface with micro-scale details, ideally suited to, for example, running simulations of droplet deposition on plant leaves.  
The surface is continuously differentiable and captures the important features of the leaf surface across multiple length scales.
The method is robust to noise, relatively performant, and most importantly avoids the storage of a whole scanned leaf at the micro-scale. 

A potential drawback of the method is its generality to other plant species. 
Wheat lends itself quite well to this method due to the almost periodic nature of the ridges on the surface. 
Other plant leaves, however, may not share this property. 
The partition of unity method is guaranteed to produce a continuous surface, though care should be taken that the blending between micro-scale segments does not introduce micro-scale features that differ too greatly from those actually present on the leaf.
For example, we took care in preparing the micro-scale dataset to avoid the case where a hair was at the edge of the segment, as the PUM blending may have distorted that hair in the multi-scale model.

It is worth noting that our particular choice of local surface reconstruction is not pivotal to this method. 
The partition of unity method could easily be applied to local interpolants $\mathcal{F}_i(\vec{x})$ constructed with any other method. 
We have used (and would recommend) the RBF method with polyharmonic splines for its robustness to noise and its minimal parameter tuning requirements.

\begin{figure}
    \centering
    \includegraphics[width=\linewidth]{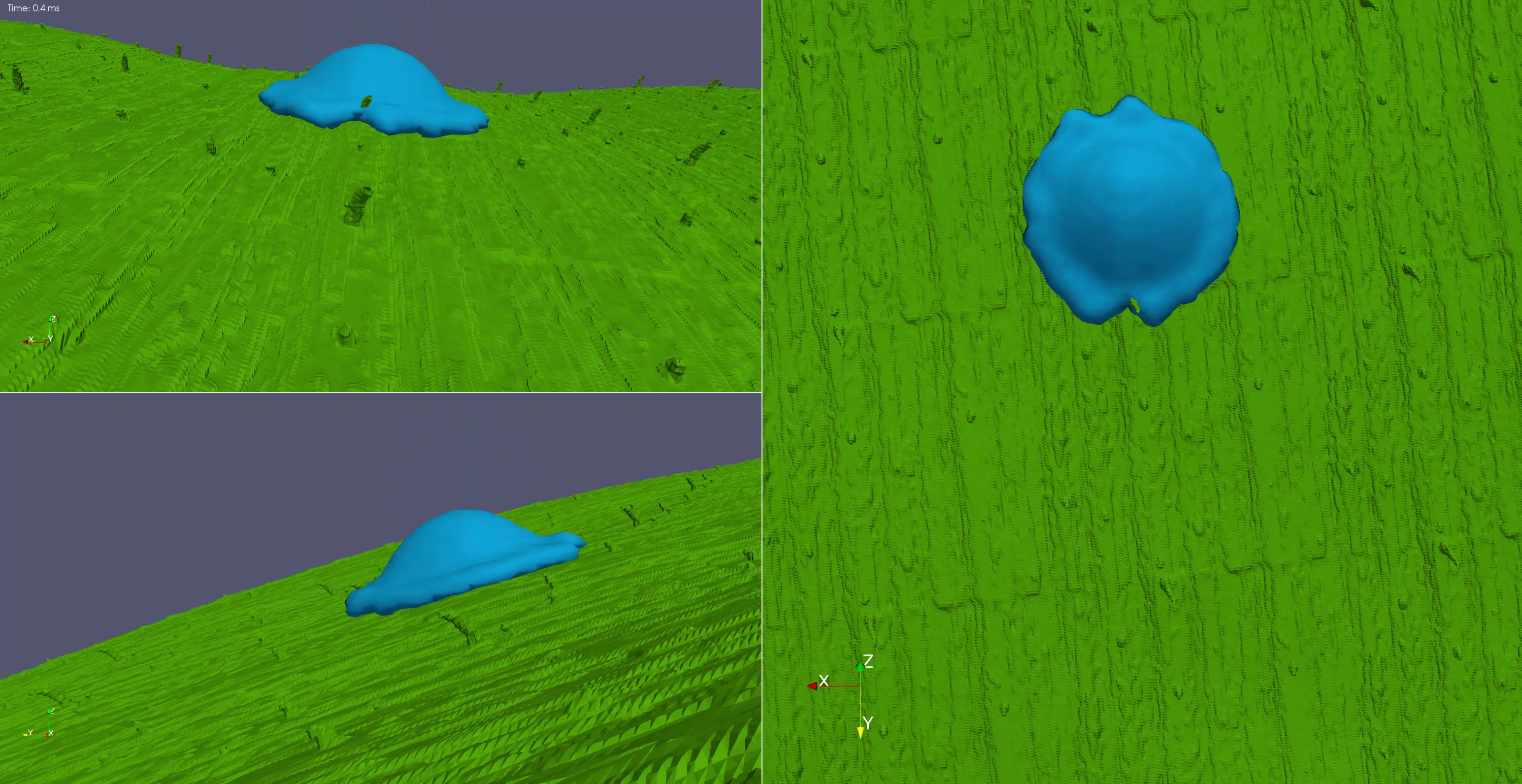}
    \caption{One frame of an impacting droplet simulation on a multi-scale wheat leaf reconstruction. The broad shape and inclination of the leaf is reconstructed from the macro-scale data, while the fine details influencing the droplet spreading are reconstructed from the micro-scale data.}
    \label{fig:dropletSim}
\end{figure}

In future work these surfaces will be used as virtual substrates on which we will study droplet deposition, dynamics, and evaporation.
Triangulations of the implicit surfaces in \cref{sec:results} could be incorporated into an existing plant retention spray model (PRSm) presented in \cite{dorr_spray_2016, PRSm2} for whole plants, or used as the domain for higher resolution studies on the individual droplet scale.
As a proof of concept, a still frame of one such simulation is shown in \cref{fig:dropletSim}. 
This simulation is generated using a smoothed particle hydrodynamics approach \cite{Monaghan_2005}, the details of which will be presented in future work.
The multi-scale nature of the virtual surface is crucial to the realism of this computational approach, permitting simulation not only of the local behaviour as the droplet spreads over the micro-scale features, but also the macro-scale movement of the droplet over the whole leaf. 

\newpage
\appendix
% \section{Nomenclature} 
% \input{src/nomenclature.tex}

\section{Defining the locally orthogonal ``leaf coordinates''}
\label{sec:orthoCoords}
Given a set of scanned data points $\vec{x}_i \in \reals^3$, we manually select control points $\vec{x}_j$ ($j = 1,\dots,n_j$) to represent the medial axis of the leaf.
We then fit a multidimensional spline, $\tilde{\vec{m}}(j)$, with the interpolation conditions $\tilde{\vec{m}}(j) = \vec{x}_j$.
The splines we use are cubic splines that are fit to the data using the julia package \texttt{Dierckx.jl} \cite{barbary_dierckxjl_2022,dierckx_curve_1995}.
We would like for the parameter of such a spline to reflect the arc length traversed between control points, so we re-parameterise by finding $t_j$ such that:
\begin{gather*}
    t_1 = 0, \\
    t_j - t_{j-1} \approx \int_{j-1}^{j} \| \tilde{\vec{m}}'(s) \|_2 \:\mathrm{d}s,
\end{gather*}
then fitting a new spline, $\vec{m}(t)$, with the interpolation conditions $\vec{m}(t_j) = \vec{x}_j$.
This spline $\vec{m}(t)$ will define the first of our leaf coordinates which indicates how far along the leaf's main axis a point lies.

The second coordinate will be normal to the leaf surface. 
For this, we fit another spline, $\tilde{\vec{V}}(t)$, through surface normals $\vec{v}_j$ at the control points.
These surface normals are estimated by principal component analysis of points $\vec{x}$ with $\|\vec{x} - \vec{x}_j\| < R$. (The constant $R$ is chosen to be roughly half the width of the leaf.)
Since both $\vec{v}_j$ and $-\vec{v}_j$ are valid surface normals, we insist that the normals agree on an orientation for the surface; that $\vec{v}_j \cdot \vec{v}_{j-1} > 0$.
We then apply the Gram-Schmidt process and define:
\begin{align*}
    \hat{\vec{V}}(t) &= \tilde{\vec{V}}(t) - \left( \tilde{\vec{V}}(t) \cdot \vec{m}'(t) \right) \vec{m}'(t),
    \\
    \vec{V}(t) &= \frac{ \hat{\vec{V}}(t) }{\| \hat{\vec{V}}(t) \|},
\end{align*}
to ensure that $\vec{V}(t) \perp \vec{m}'(t)$.
This vector will define the second leaf coordinate, indicating the height of a point with respect to the approximate plane of the leaf.

Finally, we complete the locally orthogonal coordinate system with the third leaf coordinate:
\begin{align*}
    \hat{\vec{W}}(t) &= \vec{m}'(t) \times \vec{V}(t),
    \\
    \vec{W}(t) &= \frac{\hat{\vec{W}}(t)}{\| \hat{\vec{W}}(t) \|}.
\end{align*}
We note that this coordinate system is very similar to a Frenet-Serret frame \cite{kuhnel_wolfgang_differential_2002} about the curve $\vec{m}(t)$, where $\vec{m}'(t)$, $\vec{V}(t)$, and $\vec{W}(t)$ are the tangent, normal, and binormal vectors respectively.
Our normal vector, however, is not exactly the second derivative of the curve, $\vec{m}''(t)$. 
All of these definitions culminate in a simple method to transform a point in world coordinates, $\vec{x}$, to leaf coordinates, $\vec{y}$, given in pseudocode in \cref{proc:world2leaf}.

\begin{algorithm}
    \caption{Transforming from world to leaf coordinates}
    \label{proc:world2leaf}
    \begin{algorithmic}
        \Require input point $\vec{x}$, leaf axis spline $\vec{m}(t)$, normal spline $\vec{V}(t)$, binormal spline $\vec{W}(t)$
        \Ensure leaf coordinates $\vec{y}$
        \State $t^* \gets \argmin_t \| \vec{x} - \vec{m}(t) \|^2, \quad t \in [0,t_{n_j}] $
        \State $y_1 \gets t^*$
        \State $y_2 \gets \left(\vec{x}-\vec{m}(t^*)\right) \cdot \vec{V}(t)$
        \State $y_3 \gets \left(\vec{x}-\vec{m}(t^*)\right) \cdot \vec{W}(t)$
        \State \Return $\vec{y} = [y_1, y_2, y_3]^T$
    \end{algorithmic}
\end{algorithm}

In practice, we solve the bounded minimisation problem for $t^*$ via Brent's method \cite{brent_algorithms_2013}, included in the package \texttt{Optim.jl} \cite{k_mogensen_optim_2018}.
In our experience, this does not require many iterations and the computational cost of this step is not significant compared to fitting and evaluating the radial basis function interpolants. 
The procedure to transform back, from leaf coordinate $\vec{y}$ to world coordinate $\vec{x}$, is also straightforward and efficient:
\[
    \vec{x} = \vec{m}(y_1) + y_2 \vec{V}(y_1) + y_3 \vec{W}(y_1).
\]

The notation for the leaf coordinate system is illustrated in \cref{fig:orthoCoords}.
We should be clear about the properties of the transforms between world and leaf coordinates.
If we call the world coordinate space $\mathbb{W}$ and the leaf coordinate space $\mathbb{L}$, with transform $\mathcal{T}: \mathbb{W} \rightarrow \mathbb{L}$, then the following is always true for $\vec{x} \in \mathbb{W}$:
\[
\vec{x} = \mathcal{T}^{-1}( \mathcal{T}( \vec{x} ) ),
\]
but, for $\vec{y} \in \mathbb{L}$, we cannot guarantee:
\[
\vec{y} \overset{?}{=} \mathcal{T}( \mathcal{T}^{-1}( \vec{y} ) ).
\]
This is because the transform $\mathcal{T}(\vec{x})$ involves finding the closest point $\vec{x}^* = \vec{m}(t^*)$ on the medial axis spline.
For any given point $y \in \mathbb{L}$, we do not necessarily have that $t^* = y_1$. 
This is more likely if $y_2$ or $y_3$ is large; if the point is far from the medial axis of the leaf. 
Despite this, in practice, points close together in world coordinates are close together in leaf coordinates.
Therefore, the transform is still suitable for our purposes.

\section*{Acknowledgments}
The authors gratefully acknowledge the support of the Australian Research Council through the Australian Research Council Linkage Project (LP160100707) and associated industry partners Syngenta, Nufarm, and Plant Protection Chemistry NZ Ltd.
This research was also supported by an Australian Government Research Training Program Scholarship.
The laser scan dataset of the whole wheat leaf is thanks to Dr Toby Waine and Ian Truckell at Cranfield University, School of Water, Energy and Environment.
We also acknowledge fruitful discussions with Dr Sinduja Suresh (School of Mech., Medical \& Process Engineering, Queensland University of Technology) that helped to shape the plant scanning aspects of the work.
The authors also acknowledge the support for this project provided by the Australian Research Council via the ARC Training Centre for Multiscale 3D Imaging, Modelling and Manufacturing (M3D Innovation, project IC 180100008).

\bibliographystyle{plain}
\bibliography{references.bib}

\end{document}